\documentclass[pdflatex, authoryear]{sn-jnl}

\usepackage{graphicx}%
\usepackage{multirow}%
\usepackage{amsmath,amssymb,amsfonts}%
\usepackage{amsthm}%
\usepackage{mathrsfs}%
\usepackage[title]{appendix}%
\usepackage{xcolor}%
\usepackage{textcomp}%
\usepackage{manyfoot}%
\usepackage{booktabs}%
\usepackage{algorithm}%
\usepackage{algorithmicx}%
\usepackage{algpseudocode}%
\usepackage{listings}%

\usepackage{marvosym}
\usepackage[notext]{stix2}

\theoremstyle{thmstyleone}%
\newtheorem{theorem}{Theorem}
\newtheorem{proposition}{Proposition}
\newtheorem{lemma}{\bfseries Lemma}%

\theoremstyle{thmstyletwo}%
\newtheorem{example}{\bfseries Example}%
\newtheorem{remark}{\bfseries Remark}%
\newtheorem{note}{\bfseries Note}%

\theoremstyle{thmstylethree}%
\newtheorem{definition}{Definition}%

\newcommand{\im}{\sqrt{-1}\,}
\newcommand{\ad}{{\rm ad}\,}

\newcommand{\bC}{\mathbb{C}}  
\newcommand{\bR}{\mathbb{R}}  
  
\newcommand{\bH}{\mathbb{H}}

\newcommand{\ga}{\mathfrak{a}}  
\newcommand{\gb}{\mathfrak{b}}  
\newcommand{\gc}{\mathfrak{c}}  
  
\newcommand{\gf}{\mathfrak{f}}  
\renewcommand{\gg}{\mathfrak{g}}  
\newcommand{\gh}{\mathfrak{h}}  
\newcommand{\gk}{\mathfrak{k}}  
\newcommand{\gl}{\mathfrak{l}}  
\newcommand{\gm}{\mathfrak{m}}  
\newcommand{\gn}{\mathfrak{n}}

\newcommand{\gr}{\mathfrak{r}}  
\newcommand{\gs}{\mathfrak{s}}  
\newcommand{\gt}{\mathfrak{t}}  
\newcommand{\gu}{\mathfrak{u}}

\newcommand{\gso}{\mathfrak{so}}  
\newcommand{\gaff}{\mathfrak{aff}} 
  
\newcommand{\gsu}{\mathfrak{su}}  
\newcommand{\gsl}{\mathfrak{sl}}  
\newcommand{\ggl}{\mathfrak{gl}}

\newcommand{\bt}{\begin{theorem}}  
\newcommand{\et}{\end{theorem}}  
\newcommand{\bp}{\begin{proposition}}  
\newcommand{\ep}{\end{proposition}}  
\newcommand{\bc}{\begin{corollary}}  
\newcommand{\ec}{\end{corollary}}  
\newcommand{\bl}{\begin{lemma}}  
\newcommand{\el}{\end{lemma}}  
\newcommand{\bd}{\begin{definition}}  
\newcommand{\ed}{\end{definition}} 

\newcommand{\be}{\begin{example}}  
\newcommand{\ee}{\end{example}}  
\newcommand{\br}{\begin{remark}}  
\newcommand{\er}{\end{remark}} 
\newcommand{\bn}{\begin{note}}  
\newcommand{\en}{\end{note}} 

\begin{document}

\title[Cartan Flat Non-degenerate CR Lie Groups]{Cartan Flat Non-degenerate CR Lie Groups}






\author{Keizo Hasegawa\footnote{
K. Hasegawa: Department of Mathematics, Graduate School of Science,
the University of Osaka, Toyonaka, Osaka 560-0043, Japan,
\Letter \,hasegawa@math.sci.osaka-u.ac.jp;
Faculty of Education, Niigata University, Ikarashi-nino-cho, Nishi-ku, 
Niigata 950-2181, Japan \\
}
\,and Hisashi Kasuya\footnote{
H. Kasuya:
Graduate School of Mathematics, Nagoya University,
Furo-cho, Chikusa-ku, Nagoya 464-8602, Japan, \Letter \,kasuya@math.nagaya-u.ac.jp}}

\abstract{In this paper we determine all the simply connected non-degenerate CR Lie groups, which are flat
with respect to the Cartan connection: in terms of associated Lie algebras, we assert that the only Cartan flat
non-degenerate CR Lie algebras are $\gs\gu(\textnormal{2}), \gs\gl(\textnormal{2},\bR), \ga\gf\gf (\bR) \oplus \bR$, 
and $\gh_{\textnormal{2m+1}}$
with its modifications, where $\ga\gf\gf (\bR)$ is the affine Lie algebra of dimension 2 and $\gh_{\textnormal{2m+1}}$ is
the Heisenberg Lie algebra of dimension 2m+1. Furthermore, we determine all the (flat and non-flat)
non-degenerate CR structures on each of these Lie groups.
}




\maketitle


\section{Introduction}\label{sec1}

A Sasaki manifold $M$ has a canonical underlining CR structure on its contact distribution.
Conversely, a Sasaki manifold $M$ can be characterized as a strictly pseudo-convex CR
manifold which satisfies the integrability condition (normality) to the CR structure; namely,
its canonically extended complex structure on $\tilde{M} =\bR_+ \times M$ (K\"ahler cone) is integrable,
defining a K\"ahler structure on $\tilde{M}$.

The standard Sasaki structure on $S^3 \simeq  {\rm SU}(2)$ is a (left-invariant) Sasaki structure
on ${\rm SU}(2)$, and it admits a real one-parameter family of (invariant) CR structures (obtained as a deformation)
containing the standard one but all the others are non-Sasakian. It is known (due to Rossi \cite{R}) that
these non-standard CR structures can not be globally embedded into $\bC^2$ as a hypersurface (see also Nirenberg \cite{N}).

There are many examples of Sasakian Lie algebras, but the only Sasakian unimodular Lie algebras are
$\gsu(2), \gsl(2, \bR)$ and $\gh_{2m+1}$ (Heisenberg Lie algebra) with its modifications $\overline{\gh}_{2m+1}$
(due to the work of Cort\'es and the first author \cite{CH}). In particular, the only Sasakian semi-simple Lie algebras are 
$\gsu(2), \gsl(2, \bR)$ (already known for contact structures due to Boothby-Wang \cite{BW}), 
and the only Sasakian nilpotent Lie algebras are Heisenberg Lie algebras (due to Andrada-Fino-Vezzoni \cite{AFV}).
We see in this paper (see Theorem~1 in Section~4), extending this result, that the only strictly pseudo-convex CR
nilpotent Lie algebras are Heisenberg Lie algebras.
It is known (due to Hano \cite{H}) that a K\"ahlerian  unimodular solvable Lie algebras is meta-abelian. 
To be more precise, it is a modification $\overline{\ga}_{2m}$ of an abelian Lie algebra $\ga_{2m}$,
and the Sasakian solvable Lie algebra $\overline{\gh}_{2m+1}$ above is a central extension of 
$\overline{\ga}_{2m}$ ({\cite{AHK, CH}). 


A Cartan connection for CR structures is the Cartan connection of type $G/H_0$, where
$G={\rm SU}(p+1, q+1)$ with $p+q = m$ and $H_0$ is the isotropic subgroup of $G$
for the canonical action of $G$ on the standard hyperquadric surface of $\bC^{m+1}$ 
(where $(p,q)=(m, 0)$ correspond to the sphere $S^{2m+1}$). Note that $G$ is a semisimple Lie group, and $H_0$ includes
a Borel subgroup of $G$ (Parabolic Cartan connection).
A hypersurface in $\bC^{m+1}$ has a canonical CR structure induced from the complex structure of $\bC^{m+1}$.

The theory of Cartan connections (introduced by Cartan \cite{C} and developed by Tanaka \cite{T} and 
Chern-Moser \cite{CM}) can be applied to CR structures on hypersurfaces to grasp their local invariants. In particular,
a hypersurface with a flat Cartan connection is locally isomorphic to a hyperquadric surface (which includes
a sphere). Cartan himself, applying the so-called {\em Cartan's equivalence method}, classified all the homogeneous
CR manifolds of dimension three, in particular three-dimensional CR Lie groups up to equivalence \cite{C}.
The arguments in \cite{C} on the classification were elaborated with modern terminologies in a recent paper of Bor-Jocobowitz \cite{BJ}
(also in an earlier paper of Burns-Shnider \cite{BS}).

We show in this paper a partial generalization of Cartan's result (see Theorem~4 in Section~6). 
Namely, we prove that the only Cartan flat non-degenerate CR Lie algebras are
$\gsu(2), \gsl(2, \bR), \gaff (\bR) \oplus \bR$ and $\gh_{2m+1}$ with its modifications $\overline{\gh}_{2m+1}$.
Furthermore, we determine all the (flat and non-flat) non-degenerate CR structures on each of these 
Lie groups with explicit holomorphic embeddings of hypersurface type in $\bC^{m+1}$.


\section{Preliminaries}\label{sec2}

We recall some basic terminologies and definitions for our study of CR Lie groups (Lie algebras)
in this paper.

\bd
An {\em almost Hermitian structure} on a Lie algebra $\gg$ (of even dimension)  is a pair $(\langle \cdot , \cdot \rangle, J)$ consisting of a 
scalar product $\langle \cdot , \cdot \rangle$ on $\gg$ with a skew-symmetric complex structure $J\in \mathrm{\gso}(\gg)$. The triple
$(\gg , \langle \cdot , \cdot \rangle, J)$ is an {\em almost Hermitian Lie algebra}. The almost Hermitian structure 
is {\em integrable} if the Nijenhuis tensor $N_J$ of $J$ 
vanishes, that is,
$$ N_J (X,Y) := [JX,JY] -[X,Y] -J[X,JY]-J[JX,Y] =0$$
for all $X,Y\in \gg$.

A Hermitian Lie algebra $(\gg, \langle \cdot , \cdot \rangle, J)$ is {\em K\"ahler} if its 
fundamental $2$-form $\omega =  \langle J \cdot , \cdot \rangle$ is closed. It is {\em locally 
conformally K\"ahler (shortly lcK)} if 
$$ d\omega =  \omega \wedge \theta$$ 
for some closed $1$-form $\theta \in \gg^*$ (called the {\em Lee form}). An lcK Lie algebra is {\em Vaisman}
if $\nabla \theta =0$ (parallel), where $\nabla \in \gg^* \otimes \gso (\gg)$ denotes the Levi-Civita connection. 
\ed

\bd
A {\em contact metric structure} on a Lie algebra $\gg$ of dimension $2m+1$ is
a quadrple $(\phi, \eta, \langle \cdot , \cdot \rangle, \widetilde{J})$
consisting of (i) {\em a contact structure} $\phi \in \gg^*$: $\phi \wedge (d \phi)^m \not= 0$;\;
(ii) the {\em Reeb field} $\eta \in \gg$: $i(\eta) \phi = 1, i(\eta) d \phi = 0$;
(iii) the {\em $(1, 1)$-tensor} $\widetilde{J} \, , \widetilde{J}^2 = -I + \phi \otimes \eta$;
and (iv) the {\em scalar product} $\langle \cdot , \cdot \rangle$\,: $\langle X, Y \rangle = \phi(X) \phi(Y) + d \, \phi (\widetilde{J} X, Y)$
for all $X, Y \in \gg$.

A {\em Sasaki structure} on a Lie algebra $\gg$ is
a contact metric structure
$(\phi, \eta, \langle \cdot , \cdot \rangle, \widetilde{J})$ satisfying $\langle [\eta, X], Y\rangle + \langle X, [\eta, Y]\rangle= 0$ for
all $X, Y \in \gg$ ({\em Killing field}) and such that
$J = \widetilde{J}|{\mathcal D}$ is integrable on ${\mathcal D} = \ker \phi$ ({\em CR-structure}).

For any simply connected Sasakian Lie group $G$, its {\em K\"ahler cone} $C(G)$ is defined as
$C(G)=\bR_+ \times G$ with the K\"ahler form $\Omega=r d r \wedge \phi + \frac{r^2}{2}d \phi$,
where a compatible complex structure $\widehat{J}$ is defined by 
$\widehat{J} \eta = \frac{1}{r} \partial_r$ and $\widehat{J} |{\mathcal D} = J$.


For any Sasakian Lie group $G$ with contact form $\phi$,  we can define
an lcK form $\omega=\frac{2}{r^2} \Omega = \frac{2}{r} d r \wedge \phi + d \phi$;
or taking $t= -2 \log r$,
$\omega= - d t \wedge \phi + d \phi$ on $\bR \times G$,
which is of Vaisman type. 
\ed

\bn
We can define a family of complex structures $\widehat{J}$ compatible with
$\omega$ by
$$\widehat{J} (\partial_t-l \eta)=k \eta,\; \widehat{J} (k \eta)=-(\partial_t-l \eta),\; \widehat{J} |{ D} = J, \eqno{(2.1)}$$
where $k, l \in \bR \,(k \not=0)$.
\en



\bd
A {\em CR structure} on a Lie algebra $\gg$ is a complex subalgebra $\gk$ of
$\gg_\bC$ such that $\gk \cap \overline{\gk} = \{0\}$, where  $r= \dim \gg_\bC - 2 \dim \gk \ge 0$ is
the codimension of the CR structure. Note that the case $r=0$ corresponds to the integrable 
complex structure on $\gg$.

The {\em Levi form} is a Hermitian form on $\gk$ with values in $\gg_\bC/\gk \oplus \bar{\gk}$, defined by
$h(X, Y) = - \sqrt{-1} \, [X, \bar{Y}]$ (mod $\gk \oplus \bar{\gk})$. The CR structure is {\em non-degenerate} if the Levi form
is non-degenerate, and {\em strictly pseudoconvex} if the Levi form is either positive definite or negative definite.

Two CR Lie algebras $\gg_1$ and $\gg_2$ are {\em CR-equivalent} if there is a linear isomorphism 
$\Phi: \gg_1 \rightarrow \gg_2$ such that the extended linear isomorphism $\widehat{\Phi}: (\gg_1)_\bC \rightarrow (\gg_2)_\bC$ 
maps $\gk_1$ to $\gk_2$; or equivalently, $\Phi$ maps ${\mathcal D}_1$ to 
${\mathcal D}_2$ such that $J \,\Phi_{{\mathcal D}_1} = \Phi_{{\mathcal D}_2} \,J$. Note that $\Phi$ is 
not necessarily a Lie algebra isomorphism.
\ed

{\em In this paper we consider only non-degenerate CR structures with $r=1$.}

\bn \hspace{2em}

\begin{list}{}{\topsep=3pt \leftmargin=7pt \itemindent=3pt \parsep=0pt \itemsep=5pt}

\item[(1)]
A CR structure on a Lie algebra $\gg$ is also defined by a pair $({\mathcal D}, J)$,
where $\mathcal D$ is a subspace of codimension 1 of $\gg$, and $J: {\mathcal D} \rightarrow {\mathcal D}$ with $J^2=-1$
satisfying the conditions: for $X, Y \in {\mathcal D}$, $[JX, Y] + [X, JY] \in {\mathcal D}$ and $N_J (X,Y)=0$.

\item[(2)]
We have the correspondence: $\gk = \{X - \sqrt{-1} J X | X \in {\mathcal D}\} \Leftrightarrow {\mathcal D}_\bC = \gk \oplus \overline{\gk}$.

\item[(3)]
For a strictly pseudoconvex CR structure $({\mathcal D}, J)$ on $\gg$, we have a contact metric structure 
 $(\phi, \eta, \langle \cdot , \cdot \rangle, \widetilde{J})$ with ${\mathcal D} = \ker \phi$. It defines a Sasaki structure on $\gg$
if and only if the condition $[\eta, \gk] \subset \gk$ ({\em normality}) holds. This is equivalent to say that  
the extended complex structure $\widehat{J}$ on the {K\"ahler cone} $C(G)$ is integrable.

\item[(4)]
Let $\gg$ be a Lie algebra of dimension~3. Then, CR structures $\gk$ on $\gg_\bC$ are of dimension~1, and
the set of non-degenerate CR structures $\gk$ on $\gg_\bC$ can be identified with
the elements of ${\rm P}(\gg_\bC) \cong \bC {\rm P}^2$ which are {\em regular}, meaning that $\gk \cap \overline{\gk} = \{0\}$
(non-real) and $\gk \oplus \overline{\gk}$ is not a Lie subalgebra of $\gg_\bC$ (non-degenerate). 
\end{list}

\en

\be
A {\em non-degenerate real hyperquadric} $Q$ of dimension $2m+1$ is a hypersurface in $\bC^{m+1}$ 
defined by the following equation:
$$ v = \sum_{i,j =1}^{m} h_{i,j} z_i \bar{z_j},$$
in the complex coordinate $(z_1, z_2,..., z_m, z_{m+1}),\, z_{m+1}= u + \sqrt{-1} v$ of $\bC^{m+1}$, where $H= (h_{i,j})$
is a non-degenerate Hermitian matrix.
A hyperquadric $Q$ is of {\em type $(p,q)$} $(p+q=m)$  when $H$ has $p$ positive and $q$ negative eigenvalues. 
In the complex homogeneous coordinates $(\zeta_0, \zeta_1, \zeta_2,..., \zeta_{m+1})$ of $\bC P^{m+1}$,
$\bar{Q}$ is expressed as
$$\frac{\sqrt{-1}}{2} (\zeta_0 \bar{\zeta}_{m+1} - \bar{\zeta_0} \zeta_{m+1})   =\sum_{i,j =1}^{m} h_{i,j} \zeta_i \bar{\zeta_j}, \eqno{(2.2)}$$
which is a homogeneous space of $PSU(p+1, q+1)$ with isotropy subgroup $H_0$.
\ee


\bn
The sphere $S^{2m+1}$ is a 
hyperquadric surface of type $(m, 0)$, which admits a canonical strictly pseudoconvex CR structure.
As we  see later, the Heisenberg Lie group $H_{2m+1}$ admits a non-degenerate CR structure of any type
$(p,q)$.
\en

\section{Complex structures on $\gu(\textnormal{2}), \ggl(\textnormal{2},\bR)$ and $\bR \oplus \gh_\textnormal{2m+1}$}\label{sec3}

In this section we briefly review the results of our previous work \cite{CH} on
complex structures on $\ggl(2,\bR), \gu(2)$ and $\bR \oplus \gh_{2m+1}$. We refer to the original paper
for a detailed discussion and proofs.

Let $\gg_1=\gu(2)=\bR \oplus {\mathfrak su}(2)$ and $\gg_2=\gg\gl(2,\bR)=\bR \oplus {\mathfrak sl}(2, \bR)$.
For  $\gg_1$, take a basis $\{X_1, Y_1, Z_1 \}$ for ${\mathfrak su}(2)$
with bracket multiplication defined by
$$[X_1,Y_1]=-Z_1,\, [Z_1, X_1]=-Y_1,\, [Z_1, Y_1]=X_1, \eqno{(3.1)}$$
and $T_1$ as a generator of the center $\bR$ of $\gg_1$, where we set
$$ 
X_1=\frac{1}{2}\left(
\begin{array}{cc}
0 & -1\\
1 & 0
\end{array}
\right),
\;
Y_1=\frac{1}{2}\left(
\begin{array}{cc}
0 & \sqrt{-1}\\
\sqrt{-1} & 0
\end{array}
\right),
\;
Z_1=\frac{1}{2}\left(
\begin{array}{cc}
\sqrt{-1} & 0\\
0 & -\sqrt{-1}
\end{array}
\right).
$$

For  $\gg_2$, take a basis $\{X_2, Y_2, Z_2 \}$ for ${\mathfrak sl}(2, \bR)$
with bracket multiplication defined by
$$[X_2,Y_2]=Z_2,\;[Z_2,X_2]=-Y_2,\, [Z_2, Y_2]=X_2, \eqno{(3.2)}$$
and $T_2$ as a generator of the center $\bR$ of $\gg_2$, where we set
$$ 
X_2=\frac{1}{2}\left(
\begin{array}{cc}
0 & 1\\
1 & 0
\end{array}
\right),
\;
Y_2=\frac{1}{2}\left(
\begin{array}{cc}
1 & 0\\
0 & -1
\end{array}
\right),
\;
Z_2=\frac{1}{2}\left(
\begin{array}{cc}
0 & -1\\
1 & 0
\end{array}
\right).
$$

\bl
Each $\gg_i,\, i=1,2$ admits a family of complex structures $J_{i, \delta}, \delta=k~+~\im l$ $(k, l \in \bR, k \not=0)$
defined by
$$J_{i,\delta} (T_i-l Z_i)=k Z,\; J_{i,\delta} (k Z_i)=-(T_i-l Z_i),\; J_{i,\delta} X_i= \pm Y_i,
\; J_{i,\delta} Y_i= \mp X_i.$$
Conversely, the above family of complex structures exhausts all
complex structures on each Lie algebra $\gg_i$.
\el


Let $\gg\gh=\bR \oplus \mathfrak{h}_{2m+1}$ be a nilpotent Lie algebra with 
a standard basis $\{T, X_i,Y_i,Z\}$
for which non-zero bracket multiplications are given by
$$[X_i,Y_i]=-Z,\; i=1, 2, ..., m. \eqno{(3.3)}$$

\bl
The Lie algebra 
$\gg\gh = \mathbb{R}\oplus \mathfrak{h}_{2m+1}$ 
admits a family of complex structures 
$J=J_{(\delta;\, \varepsilon_1, \varepsilon_2,..., \varepsilon_m)}, \delta=k+\im l\in \bC\;(k, l \in \bR, k \not=0), \varepsilon_i=\pm 1$, 
defined by
$$J (T-l Z)=k Z,\; J (k Z)=-(T-l Z),\; J X_i= \varepsilon_i Y_i,\; J Y_i= -\varepsilon_i X_i.$$
Conversely, the above family of complex structures exhausts all
complex structures on $\gg\gh$.
\el
\bp
(\cite{CH}) The complex structures on $\gg\gl(2, \bR), \gu(2)$ and $\gg\gh_{2m+2}$ 
are biholomorphic to $\bC \times \bH, \bC^2 \backslash \{0\}$ and $\bC^{m+1}$ respectively.
\ep

\begin{list}{}{\topsep=3pt \leftmargin=7pt \itemindent=3pt \parsep=0pt \itemsep=5pt}

\item[(1)] 
For the case of $\bR \times \widetilde{SL}(2, \bR)$,
we consider a biholomorphic map $\Phi_\delta$:
$$\Phi_\delta: (\bR \times {SL}(2, \bR), J_\delta) \longrightarrow \bH \times \bC^*$$
defined by
$$ g = (t, \left(
\begin{array}{cc}
a & b\\
c & d
\end{array} 
\right))
\longrightarrow ( \left( \frac{ac +bd}{c^2+d^2} \right) + \sqrt{-1} \left( \frac{1}{c^2+d^2} \right),\; e^{\delta\, t} \,(d + \sqrt{-1}\, c) ), $$
where $\delta=k+\sqrt{-1}\, l$. 
Then, it induces a biholomorphic map $\bar{\Phi}_\delta$
from the universal coverings $(\bR \times \widetilde{SL}(2, \bR), \bar{J}_\delta)$
to $\bH \times \bC$, where $\bar{J}_\delta$ is the induced integrable complex structure from $J_\delta$.

\item[(2)]
 For the case of $\bR \times SU(2)$,
we have a canonical biholomorphic map $\Phi_\delta$:

$$\Phi_\delta: (\bR \times SU(2), J_\delta) \longrightarrow \bC^2 \backslash \{0\}$$
defined 
$$(t,z_1,z_2) \longrightarrow (e^{\delta \,t} z_1, e^{\delta \,t} z_2),$$

\noindent where $\delta=k+{\im l}$ and
$SU(2)$ is identified with
$$S^3=\{(z_1,z_2) \in \bC^2 \,|\;|z_1|^2+|z_2|^2=1\}$$
by the correspondence
$$
\left(
\begin{array}{cc}
z_1 & -\overline{z}_2\\
z_2 & \overline{z}_1
\end{array}
\right)
\longleftrightarrow (z_1,z_2).
$$

\item[(3)]
For the case of $GH_{2m+2}=\bR \times H_{2m+1}$,
we have a canonical biholomorphic map $\Phi_\delta$:
$$\Phi_\delta: (\bR \times H_{2 m+ 1}, J_\delta) \rightarrow \bC^{m+1}$$
defined by

$$\big(t,\, \left(
\begin{array}[c]{ccc}
1 & {\bf x} & z\\
0 & {\rm I}_m & {\bf y}^t\\
0 & 0 & 1
\end{array}
\right) \,\big)
\rightarrow
({\bf x} + \sqrt{-1} {\bf y},\, (2 kt + \frac{1}{2}(\|{\bf x}\|^2 + \|{\bf y}\|^2)) + \sqrt{-1}\,
(2 (lt+z) - {\bf x} \cdot {\bf y})),
$$
where ${\bf x}=(x_1, x_2, ..., x_m), \, {\bf y}=(y_1, y_2, ..., y_m) \in \bR^m$, ${\bf x} \cdot {\bf y} =
\sum_{i=1}^{m} \varepsilon_i \, x_i y_i$ and $\|{\bf x}\|^2 = {\bf x} \cdot {\bf x}$.

\end{list}

\section{CR structures on $\gs\gu(\textnormal{2})$, 
$\gs\gl(\textnormal{2},\bR), \ga \gf \gf (\bR) \oplus \bR$ and $\gh_\textnormal{2m+1}$}\label{sec4}

Non-degenerate CR structures on $\gs\gu(2)$, $\gs\gl(2,\bR), \ga \gf \gf (\bR) \oplus \bR$ are
known \cite{C, BJ}. We determine in this section non-degenerate CR structures on $\gh_{2m+1}$, and prove a structure theorem 
of non-degenerate CR nilpotent Lie algebra $\gg$ in general (Theorem~1).

\begin{list}{}{\topsep=3pt \leftmargin=7pt \itemindent=3pt \parsep=0pt \itemsep=5pt}

\item[(1)] 
For $\gg={\mathfrak su}(2)$, the Cartan-Killing form is negative definite.
Since we have $\mathrm{Aut} (\gg) \cong SO(3)$, we can take any subspace of dimension 2
for ${\mathcal D}$. So, we can let ${\mathcal D} = <X, Y>_\bR$ with a family of complex structures 
defined by $J X = tY \;(t > 0)$. Then,
each $\gk_t = <X + t \sqrt{-1} Y>_\bC$ defines a strictly pseudoconvex CR structure on $\gg$. Since we have
$$[X + t \sqrt{-1} Y, Z] = Y  - t \sqrt{-1} X = -t \sqrt{-1} \,(X + \frac{\sqrt{-1}}{t}Y),$$
which is in $\gk_t$ iff $t=1$. That is, the normality is satisfied only for $t=1$, which is the case of the standard Sasaki structure on $S^3$.
Note that it is well known that CR structure $\gk_t, t \not=1$  on $S^3$ is non-embeddable in $\bC^2$ as a hypersurface.

\item[(2)] 
For $\gg ={\mathfrak {sl}}(2,\bR)$, the Cartan-Killing form is Lorentzian. Since we have $\mathrm{Aut} (\gg) \cong SO(2,1)$,
we can take, as a subspace of dimension 2 for CR structure on $\gg$, either ${\mathcal D}_Z = <X, Y>_\bR$ 
or  \mbox{${\mathcal D}_X = <Y, Z>_\bR$},
where $Z$ is a timelike vector, and $X, Y$ are spacelike vectors with respect to the Lorentzian metric.
Then, each \mbox{$\gk_{Z, t} = <X - t \sqrt{-1} Y>_\bC$} or \mbox{$\gk_{X, t} = <Y - t \sqrt{-1} Z>_\bC$} with $t > 0$ defines
a strictly pseudoconvex CR structure on $\gg$. The normality is satisfied only for $\gk_{Z, 1} = <X -  \sqrt{-1} Y>_\bC$,
which is the case of the standard Sasaki structure on ${\textrm SL}(2,\bR)$.

\item[(3)] For $\gg= \ga \gf \gf (\bR) \oplus \bR$, take a basis a basis $\{X, Y, Z\}$ for which the bracket multiplication is
defined by
$$[X,Y]=-Z,\; [Y, Z]=0,\; [Z,X]=Z. \eqno{(4.1)}$$
Take another basis, changing $Z$ into $Z'=Z-Y$, we have
$$[X,Y]=-Z'-Y,\; [Y, Z']=0, \; [Z',X]=0.$$
Let $x, y, z'$ be the Maurer-Cartan forms corresponding to $X, Y, Z'$ respectively. We have
$$d x = 0, \; d y = x \wedge y,\; d z' = x \wedge y.$$
Then, we have only the CR structure ${\mathcal D} = <X, Y>_\bR$ and
\mbox{$\gk = <X -  \sqrt{-1} Y>_\bC$} with the contact form $\phi = z'$.
Note that this CR structures is equivalent to the CR structure $\gk_{Z, 1}$ on
${\mathfrak {sl}}(2,\bR)$ by the identity linear map.

\item[(4)] For $\gg=\gh_{2m+1}$, since the CR structure is non-degenerate and $Z(\gg) = <Z>_\bR$,
we see that $Z(\gg) \cap  {\mathcal D} = \{0\}$ and $Z$ must be the Reeb field of the contact structure
${\mathcal D} = <X_i, Y_j>~, i, j=1,2,.., m$.
Then the normality is trivially satisfied, extending any complex structure on ${\mathcal D}$ to 
$\gg\gh = \bR \oplus \gh_{2m+1}$. According to Lemma 2, the complex structures on ${\mathcal D}$
are given by $JX_i= \varepsilon_i  Y_i,\; i=1,2,..,m,\; \varepsilon_i = \pm 1$. 
We see that the above biholomorphism $\Phi_1$ from $GH_{2m+2}$ to $\bC^{m+1}$ embeds 
the CR structure on the Heisenberg Lie group $H_{2m+1}$ into a hyperquadratic surfaces of all type $(p, q)$ in $\bC^{m+1}$. 
And conversely, these exhausts all non-degenerate CR structures on $H_{2m+1}$.

\end{list}


\bt
A non-degenerate CR nilpotent Lie algebra $\gg$ is pseudo-Sasakian,
which is a central extension by $\bR$ of a pseudo-K\"ahlerian nilpotent Lie algebra.
In particular, a strictly pseudo-convex CR nilpotent Lie algebra $\gg$ is Sasakian, and thus a Heisenberg Lie algebra
(which is a central extension by $\bR$ of an abelian Lie algebra).
\et

\begin{proof}
Let ${\mathcal D} \subset \gg$ be a non-degenerate CR structure on $\gg$ of dimension $2m+1$.
with a contact form $\phi$ (${\mathcal D} = \ker \phi$) and
the Reeb field $\eta$. We have $\gg = {\mathcal D} \,\oplus <\eta>$, and $i(\eta)\phi=1, i(\eta)d \phi =0$. 
Let $Z(\gg)$ be the center of $\gg$ (which is non-trivial). It is sufficient to show that $Z(\gg) = <\eta>$, 
since it shows that the integrability condition
is trivially satisfied, and thus $\gg$ is a pseudo-Sasakian nilpotent Lie algebra.  In fact, we then have
$$ 0 \rightarrow \bR \rightarrow \gg \rightarrow \gn \rightarrow 0,$$
where the induced non-degenerate closed $2$-form $d\overline{\phi}$ on $\gn = \gg/Z(\gg)$
defines a pseudo-K\"ahlerian structure on $\gn$, and thus $\gg$ is pseudo-Sasakian. In particular, in case $\gg$ is a 
strictly pseudo-convex CR nilpotent Lie algebra, $d\overline{\phi}$ is K\"ahlerian, and thus
$\gn$ must be abelian, that is, $\gg$ is  a Heisenberg Lie algebra.
 

Now, since $d \phi$ is non-degenerate on
$\mathcal D$ and $d \phi(X, Y) = \phi([X,Y])$, we have $Z(\gg) \cap {\mathcal D} = \{0\}$ and
$\gg = {\mathcal D} \,\oplus Z(\gg)$. Since $i(\eta)d \phi =0$, we must have $Z(\gg) = <\eta>$.
\end{proof}

\bn
There is a pseudo-Sasakian nilpotent Lie algebras as a central extension by $\bR$ of
a pseudo-K\"ahler (non-abelian) nilpotent Lie algebra. For example, we can take
$\gg\gh = \bR \oplus \gh_3$ as a (non-K\"ahlerian) pseudo-K\"ahlerian nilpotent Lie algebra of dimension $4$.
\en

\section{Cartan Connection}\label{sec5}

Let $M$ be a manifold of dimension $n$ and let $P=L(M)$ be the frame bundle over $M$ with group $H=GL(n, \bR)$ and
the projection $\pi: P \rightarrow M$.
The {\em solder form} $\theta$ on $P$ is defined by
$$\theta(X) = u^{-1}(\pi_* (X)) \in \bR^n$$ 
for $X \in T_u P \;(u \in P)$. 

Let $w$ be a linear connection on $P$ satisfying
\begin{list}{}{\topsep=5pt \leftmargin=17pt \itemindent=3pt \parsep=0pt \itemsep=3pt}
\item[(i)] $R^*_a w = a^{-1} w a \; (a \in H)$
\item[(ii)] $w(A^*) = A \; (A \in \gh)$;  $A^*$: the fundamental vector field on $P$
\end{list}

Then we have 
$$ T_u (P) \; \widetilde{\rightarrow} \; \bR^n \oplus \gh$$
by $X \rightarrow \omega(X)=(\theta(X), w(X))$ (absolute parallelism).

Let $G = {\rm Aff} (n, \bR) = \bR^n \rtimes H$ with Lie algebra $\gg = \bR^n \oplusrhrim \gh$, and 
$A(M)=P \times_H G$ (affine frame bundle over $M$).
$P=L(M)$ is a $H$-reduction of $A(M)$; and we have $T(M) \; \widetilde{\rightarrow} P \times_H \gg/\gh$
with adjoint representation of $H$ on $\gg/\gh$.
\vspace{0.7em}

We see that $\omega$ defines a {\em Cartan connection} on $P$ as defined in the following.


\bd
Let $H$ a closed subgroup of $G$ and let $P$ be a $H$-principal bundle over a manifold $M$
with a $\gg$-valued $1$-form $\omega$ ({\em Cartan connection}), satisfying 

\begin{list}{}{\topsep=5pt \leftmargin=10pt \itemindent=3pt \parsep=0pt \itemsep=5pt}
\item[(i)] $R^*_a \omega = {\rm Ad} (h^{-1}) \,\omega \; (h \in H)$
\item[(ii)] $\omega (A^*)= A \;(A \in \gh)$
\item[(ii)] $\omega_p: T_p (P) \; \widetilde{\rightarrow} \; \gg \;(p \in P)$ \;(absolute parallelism)
\end{list}

\bn
We have ${\rm dim}\, M = {\rm dim}\, G - {\rm dim}\, H$; $P$ is a $H$-reduction of the $G$-principal bundle
$P \times_H G$; and $T(M) \;\widetilde{\rightarrow} \; P \times_H \gg/\gh$. In particular, the Cartan connection 
$\omega \;{\rm mod} \; \gh: P \rightarrow \gg/\gh$ is horizontal and defines a {\em solder form} on $P$.
\en

The curvature $\Omega$ of the Cartan connection $\omega$ is defined as
$$ \Omega = d \omega + \omega \wedge \omega.$$

\ed


\be
Let $M = G/H$ be a homogeneous space. Then, $G$ is a $H$-principal bundle over $M$. The Maurer-Cartan form 
$\omega$ ($\gg$-valued $1$-form) on $G$ defines a flat Cartan connection with the structure condition:
$$d \omega = - \omega \wedge \omega.$$
\ee


\bn
In case $M=G/H$ is {\em reductive}; that is,  $\gg$  has the decomposition
$$\gg=\gh \oplus \gk$$
with $[\gk, \gh] \subset \gk$ for some closed subgroup $K$ of $G$,
the Maurer-Cartan form $\omega$ decomposes as
$$\omega = \omega_\gh + \omega_\gk,$$
where $\omega_\gh$ defines a connection on $G$ as a $H$-principal bundle, and $\omega_\gk$ defines
a solder form.
\en


\bn
A Cartan connection $\omega$ of type $G/H$ on a $H$-principal bundle $P$ is of {\em reductive type}
if $G/H$ is of reductive type. Affine connection is of reductive type.
\en

{\em Cartan's equivalence method is to associate any non-degenerate CR manifold $M$ 
(of hypersurface type) $H_0$-principal bundle $P$ over $M$ with a Cartan connection $\omega$ with value in $\gs\gu(p+1, q+1)$,
where $H_0$ is the isotropy subgroup of $SU(p+1,q+1)$ as acting on a hyperquadratic
surface of type $(p,q)$. The existence and uniqueness of the bundle $P$ is due to Cartan, Tanaka, and Chern-Moser \cite{C, T, CM}.
To be more precise, for given non-degenerate CR manifolds $M_1, M_2$ and a CR-equivalent map $\phi: M_1 \rightarrow M_2$,
there exists uniquely a bundle isomorphism $\tilde{\phi}: P_1 \rightarrow P_2$ preserving the Cartan connections
$\omega_1, \omega_2$ on the corresponding principal bundles $P_1, P_2$ respectively
(see Th. 4.4, \cite{T}).

}


\section{Cartan flat non-degenerate CR Lie groups}\label{sec6}

Let $G= SU(p+1,q+1)$ with Lie algebra $\gg=\gs\gu(p+1, q+1), p+q = m$, and  $H_0$ the isotropy
subgroup of $G$  with  Lie algebra $\gh_0$.
Since the Cartan connection $\omega$ on a CR manifold $M$  gives $T(M) \;\widetilde{\rightarrow} \;P \times_H \gg/\gh_0$; and 
$\omega \,{\rm mod} \, \gh_0: P \rightarrow \gg/\gh_0$ is horizontal, defining a solder form on $P$, 
we have a canonical isomorphism as vector spaces:
$$T_x(M) \hookrightarrow T_p(P) = \gg \rightarrow \gg/\gh_0,$$
where $x \in M, p \in P$ with $\pi(p)=x$.
In particular, if the Cartan connection is flat, the horizontal distribution is integrable.
Hence, for  a flat non-degenerate CR Lie group $L$ with Lie algebra $\gl$,
we have a Lie algebra structure on $\gg/\gh_0$ and a canonical isomorphism of Lie algebras
$$ \gl \hookrightarrow \gg \rightarrow \gg/\gh_0.$$
Here, we take an open neighborhood $U$ of $e$ on $L$, setting $\gl = T_eU$ and $\omega_0$ the restriction of $\omega$
on $\gl$. Then,
$\omega_0: \gl \hookrightarrow \gg$ is injective as we have ${\rm dim}\, \gl = {\rm dim}\, \gg/\gh_0$;
and it is a Lie algebra homomorphism as 
the ${\rm graph} (\omega_0) = \{(X, \omega_0(X)) \in \gl \oplus \gg\, | X \in \gl\}$ is a Lie subalgebra of $ \gl \oplus \gg$
due to the $L$-invariance and flatness of the connection (i.e. involutivity of  ${\rm graph} (\omega_0)$).
In fact, applying the Cartan's equivalence method, the action of $L$ as a CR automorphism on $L$ lifts to
that as a bundle automorphism on $P$ preserving the connection $\omega$, and the vector field
$(X, \omega (X))$ on $TU \times TG$ is involutive.  Hence,
${\rm graph} (\omega_0) = \{(X, \omega_0(X)) \in \gl \oplus \gg\, | X \in \gl\}$ defines a Lie subalgebra of $\gl \oplus \gg$,
showing that $\omega_0: \gl \hookrightarrow \gg$ is an injective
Lie algebra homomorphism.
\vspace{0.3em} 

We have thus shown the following fundamental result.


\bt
Let  $L$ be a non-degenerate flat CR Lie group with its Lie algebra $\gl$.
Then $\gl$ must be a subalgebra of $\gg$ satisfying
$$\gg= \gl \oplus \gh_0  \eqno(6.1)$$
for $\gg = \gs\gu(p+1, q+1)$ with its isotropy subalgebra $\gh_0$ as a vector space direct sum.
\et

Note that
$\gg$ has the decomposition
$$\gg= \gh_{2m+1} \oplus \gh_0. \eqno(6.2)$$

\bn
In terms of Lie groups, this is equivalent to say that the Lie group $L$ acts simply transitively 
on the homogeneous space $G/H_0$. 
\en



Let us consider first the case $p=1, q=0$, that is,  $\gg = \gs\gu(2, 1)$. We have

$$\gg= \{X \in \gs\gl (3, \bC) |\;
\left(
\begin{array}{ccc}
0 & 0 & \sqrt{-1}\\
0 & 1 & 0\\
-\sqrt{-1}& 0 & 0
\end{array} 
\right) X + X^* 
\left(
\begin{array}{ccc}
0 & 0 &  \sqrt{-1}\\
0 & 1 & 0\\
- \sqrt{-1}& 0 & 0
\end{array} 
\right) =0
\}
$$
By calculation, $X \in \gg$ can be written as
$$X= X_{(a, b, \alpha, \beta, \gamma)}=   \left(
\begin{array}{ccc}
\beta & \gamma & b\\
\alpha & -(\beta - \overline{\beta}) & -  \im \overline{\gamma}\\
a &  \im \overline{\alpha} & - \overline{\beta}
\end{array} 
\right), \alpha, \beta, \gamma \in \bC, a, b \in \bR.$$
The isotropy subalgebra $\gh_0$ of $\gg$ is given by
$$\gh_0=  \Bigl\{\left(
\begin{array}{ccc}
\beta & \gamma& b\\
0& -(\beta - \overline{\beta})& -  \im \overline{\gamma}\\
0 & 0 & - \overline{\beta}
\end{array} 
\right) \big|
\; \beta, \gamma \in \bC, b \in \bR \Big\}.$$
The Cartan subalgebra $\ga = \ga_{1,0}$ of $\gg$ is given by
$$\ga = \ga (\beta) = \Big\{\left(
\begin{array}{ccc}
\beta & 0 & 0\\
0 & -(\beta - \overline{\beta}) & 0\\
0 & 0 & - \overline{\beta}
\end{array} 
\right) \big|
\;\beta \in \bC \Big\}.
$$
The Borel subalgebra $\gb= \gb_{1,0} $ of $\gg$ is given as
$$\gb= \gb(a,\alpha,\beta) = \Big \{\left(
\begin{array}{ccc}
\beta & 0 & 0\\
\alpha & -(\beta - \overline{\beta}) & 0\\
a &  \im \overline{\alpha} & - \overline{\beta}
\end{array} 
\right) \big|
\; \alpha, \beta \in \bC, a \in \bR \Big\};
$$
that is, $\gb = \gh_3 \oplusrhrim_{ad} \ga $, where $\ga$
acts on $\gh_3$ as the adjoint representation $ad$ in $\gg$.



\noindent For the general cases: $\gg = \gs\gu(p+1, q+1)$, we have
$$\gg= \{X \in \gs\gl (m+2, \bC) |
\left(
\begin{array}{ccc}
0 & 0 & \sqrt{-1}\\
0 & I_{p, q} & 0\\
-\sqrt{-1}& 0 & 0
\end{array} 
\right) X + X^* 
\left(
\begin{array}{ccc}
0 & 0 &  \sqrt{-1}\\
0 & I_{p,q} & 0\\
- \sqrt{-1}& 0 & 0
\end{array} 
\right) =0
\}.
$$
By calculation,  $X \in \gg$ can be written as
$$X= X_{(a,b, \beta, u, v)} =\left(
\begin{array}{ccc}
\beta & - \im v^* & b\\
u & \gu_{p,q} & v\\
a &  \im \overline{u}^* & - \overline{\beta}
\end{array} 
\right),
\;\; \beta \in \bC, a, b \in \bR, u, v \in \bC^m, \eqno{(6.3)}$$
and $\gu_{p,q} = \gu (p,q)$ with the condition ${\rm tr}\; \gu_{p,q}+ 2 \, {\rm Im}\, \beta =0$.
The isotropy subalgebra $\gh_0$ of $\gg$ is given by
$$\gh_0=  \Bigl\{\left(
\begin{array}{ccc}
\beta & - \im v^* & b\\
0& \gu_{p,q} & v\\
0 & 0 & - \overline{\beta}
\end{array} 
\right) \big|
\; \beta \in \bC, a, b \in \bR, u, v \in \bC^m \Big\}.$$
The Cartan subalgebra $\ga = \ga_{p,q}$ of $\gg$ is given by
$$\ga_{p,q}= \ga (c, d, D) = \Big \{\left(
\begin{array}{ccc}
c - \im d & 0 & 0\\
0& D & 0\\
0 & 0 & - c - \im d
\end{array} 
\right) \big|\; c, d \in \bR, D \in  \gg\gl(m, \bC) \Big\};
$$
where $D = \im {\rm dia} (d_1, d_2,..., d_m), d_i \in \bR$ with
$d_1+d_2+ \cdots +d_m=2d$.
The Borel subalgebra $\gb_{p,q}$ of $\gs\gu(p+1, q+1)$ is given as
$$\gb_{p,q}= \gb (a,c, d, u, D) =  \Big\{ \left(
\begin{array}{ccc}
c - \im d & 0 & 0\\
u & D & 0\\
a &  \im \overline{u}^* & - c - \im d
\end{array} 
\right) \big| \;a, c, d \in \bR, u \in \bC^m \Big\};
$$
that is, $\gb_{p,q} = \gh_{2m+1} \oplusrhrim_{ad} \ga_{p,q}$, 
where $\ga_{p,q}$ acts on $\gh_{2m+1}$ as the adjoint representation $ad$ in $\gg$.
Note that we can decompose $U \in \ga_{p,q}$ as
$$U = U_{(c, d, D_0)}= \left(
\begin{array}{ccc}
c & 0 & 0\\
0& 0 & 0\\{}
0 & 0 & -c
\end{array} 
\right)
+
\left(
\begin{array}{ccc}
-\im d& 0 & 0\\
0 & \frac{2 \im d}{m} I_m & 0\\
0 & 0 & - \im d
\end{array} 
\right)
+
\left(
\begin{array}{ccc}
0& 0 & 0\\
0 & D_0 & 0\\
0 & 0 & 0
\end{array} 
\right),
$$
where $c, d \in \bR, D_0=\im {\rm dia}(d_1,d_2,...,d_m)$ with $d_1+d_2+ \cdots +d_m=0$.

Note that a subalgebra $\gl$ of  $\gb_{p,q}$, satisfying the condition (6.1) can be considered as
a {\em modification} $\gh_{2m+1} (\tau)$ of $\gh_{2m+1}$ by a Lie algebra homomorphism $\tau: \gh_{2m+1} \rightarrow \ga_{p,q}$,
where $\gh_{2m+1} (\tau)$ is the subalgebra of $\gh_{2m+1} \oplusrhrim_{ad} \ga_{p,q}$ consisting of the graph $ (X, \tau(X))$
with the bracket multiplication defined by 
$$[(X, \tau(X)), (Y, \tau(Y)]_\tau = ([X, Y] +\tau(X)(Y) - \tau(Y)(X), [\tau(X), \tau(Y)]).$$
Since  $\ga_{p,q}$ is abelian, we can consider $\gh_{2m+1} (\tau)$ as a Lie algebra on the space $\gh_{2m+1}$ with the modified
bracket multiplication defined by
$$[X, Y]_\tau = [X, Y] +\tau(X)(Y) - \tau(Y)(X).$$

\vspace{1em}

We will now determine all subalgebras $\gl$ of $\gg = \gs\gu(p+1, q+1)$ satisfying (6.1), and then see
which of these Lie algebras admit non-degenerate CR structures. 


We first show in {\em Example~3}  below that $\gs \gu(2), \gs \gl (2,\bR)$ and $\gh_{3}$ actually 
admit flat non-degenerate 
CR structures, explicitly constructing the subalgebras $\gl$ satisfying (6.1) which are isomorphic to these Lie algebras.

 
\be  \hspace{2em}

\begin{list}{}{\topsep=3pt \leftmargin=7pt \itemindent=3pt \parsep=0pt \itemsep=5pt}

\item[(1)] 
For $\gl={\mathfrak su}(2)$, ${\mathcal D} = <X, Y>_\bR$ and  $\gk_t = <X + t \sqrt{-1} Y>_\bC \;(t > 0)$ define
a strictly pseudoconvex CR structure on $\gl$. Let $x, y, z$ be the Maurer-Cartan form dual to $X,Y,Z$. Then we have
$$d x =  -y \wedge z, d y =  -z \wedge x, d z =  -x \wedge y$$
We set 
$$\phi = 2 z,\; \psi = \sqrt{t} x + \sqrt{-1} \frac{1}{\sqrt{t}}y.$$
Then ${\mathcal D} = {\rm ker}\; \phi$ and $\psi (\gk_t) = 0$, and we have
$$d \phi = \im \psi \wedge \overline{\psi},\;\;
d \psi = b_t \,\phi \wedge \psi + c_t \,\phi \wedge \overline{\psi},$$
where $b_t= {\displaystyle \frac{\sqrt{-1}}{4} (\frac{1}{t} + t)}, c_t = {\displaystyle \frac{\sqrt{-1}}{4} (\frac{1}{t} - t)}$.

We define $\gs\gu(2, 1)$-valued $1$-form $\Phi_t$ by
$$\Phi_t = \left(
\begin{array}{ccc}
b_t \,\phi &  \im (b_t \overline{\psi} + c_t \psi) & (|c_t|^2 - |b_t|^2) \,\phi\\
\psi & -2 b_t \,\phi & b_t \, \psi + c_t \,\overline{\psi}\\
\phi &  \im \overline{\psi} & b_t \,\phi
\end{array} 
\right),
$$
which defines CR structures with parameter $t > 0$ on $\gl$.
We see that $\Phi_t$ is closed by the bracket multiplication on $\gs\gu(2, 1)$ only for $t=1$,
defining a subalgebra isomorphic to $\gl$:

$$\Phi_1= \left(
\begin{array}{ccc}
\frac{\im}{2} \phi & - \frac{1}{2} \overline{\psi} & - \frac{1}{4} \phi\\
\psi & - \im \phi & \frac{\im}{2}\psi \\
\phi &  \im \overline{\psi} & \frac{\im}{2}\phi
\end{array} 
\right)
$$

\item[(2)] 
For $\gl={\mathfrak sl}(2, \bR)$, ${\mathcal D} = <X, Y>_\bR$ and  $\gk_t = <X - t \sqrt{-1} Y>_\bC \;(t > 0)$ define
a strictly pseudoconvex CR structure on $\gl$. Let $x, y, z$ be the Maurer-Cartan form dual to $X,Y,Z$. Then we have
$$d x =  y \wedge z, d y =  z \wedge x, d z =  - x \wedge y$$

We set 
$$\phi = 2 z,\; \psi = \sqrt{t} x - \sqrt{-1} \frac{1}{\sqrt{t}}y.$$
Then ${\mathcal D} = {\rm ker}\; \phi$ and $\psi (\gk_t) = 0$, and we have
$$d \phi = \im \psi \wedge \overline{\psi},\;\;
d \psi = b_t \,\phi \wedge \psi + c_t \,\phi \wedge \overline{\psi},$$
where $b_t= {\displaystyle -\frac{\sqrt{-1}}{4} (\frac{1}{t} + t)}, c_t = {\displaystyle -\frac{\sqrt{-1}}{4} (\frac{1}{t} - t)}$.

We define $\gs\gu(2, 1)$-valued $1$-form $\Phi_t$ by
$$\Phi_t = \left(
\begin{array}{ccc}
b_t \,\phi &  \im (b_t \overline{\psi} + c_t \psi) & (|c_t|^2 - |b_t|^2) \,\phi\\
\psi & -2 b_t \,\phi & b_t \, \psi + c_t \,\overline{\psi}\\
\phi &  \im \overline{\psi} & b_t \,\phi
\end{array} 
\right),
$$
which defines CR structures with parameter $t > 0$ on $\gl$.
We see that $\Phi_t$ is closed by the bracket multiplication on $\gs\gu(2, 1)$ only for $t=1$,
defining a subalgebra isomorphic to $\gl$:

$$\Phi_1= \left(
\begin{array}{ccc}
-\frac{\im}{2} \phi & \frac{1}{2} \overline{\psi} & - \frac{1}{4} \phi\\
\psi &  \im \phi & -\frac{\im}{2}\psi \\
\phi &  \im \overline{\psi} & -\frac{\im}{2}\phi
\end{array} 
\right)
$$

\item[(3)]
 For $\gl= \gh_3$ (Heisenberg Lie algebra), ${\mathcal D} = <X, Y>_\bR$ and  $\gk = <X + \sqrt{-1} Y>_\bC$ define
a strictly pseudoconvex CR structure on $\gl$. Let $x, y, z$ be the Maurer-Cartan form dual to $X,Y,Z$. Then we have
$d z =  x \wedge y, d x =  d y =0$. 
We set
$$\phi = z,\; \psi=\frac{1}{2} (x + \sqrt{-1} y).$$
Then ${\mathcal D} = {\rm ker}\; \phi$ and $\psi (\gk) = 0$, and we have
$$d \phi = \sqrt{-1} \psi \wedge \overline{\psi},\;\;
d \psi = 0$$

We define $\gs\gu(2, 1)$-valued $1$-form $\Phi$ by
$$\Phi = \left(
\begin{array}{ccc}
0 & 0 & 0\\
\psi & 0 & 0\\
\phi &  \sqrt{-1} \, \psi & 0
\end{array} 
\right),
$$
which defines a CR structure on $\gl$.

\end{list}

\ee


\bt
The only three-dimensional Lie algebras which admit Cartan flat non-degenerate CR 
structures are $\gs \gu(2), \gsl (2,\bR), \ga \gf \gf (\bR) \oplus \bR$ and $\gh_3$, among which
the CR structures on $\ga \gf \gf (\bR) \oplus \bR$ and  $\gsl (2,\bR)$ are equivalent,
and no other CR structures are equivalent to each other.
\et

\begin{proof}

We have seen in {\em Example 3} above that $\gs \gu(2), \gsl (2,\bR)$ and $\gh_3$ admit Cartan flat 
non-degenerate CR structures, expressing these as subalgebras satisfying the condition (6.1).
We have also seen in Section~4 that the only non-degenerate CR structure on  $\gaff (\bR) \oplus \bR$ is
equivalent to a non-degenerate CR structures of type ${\mathcal D}_Z$ on $\gsl (2,\bR)$.
Since $\gs \gu(2)$ and $\gsl (2,\bR)$ are the only three-dimensional semisimple Lie algebras,
it is sufficient to show that the only three-dimensional solvable Lie algebras $\gl$ satisfying the condition (6.1)
are $\gaff (\bR) \oplus \bR$ and $\gh_3$. We actually show the stronger assertion that the only three-dimensional 
non-degenerate CR solvable Lie algebra $\gl$ contained in the Borel subalgebra $\gb$ is  either $\gaff (\bR) \oplus \bR$ or $\gh_3$. 
Note that since the Borel subalgebras are all conjugate and any solvable
Lie subalgebra is contained in a Borel subalgebra, we may actually assume that
$\gl$ is a subalgebra of $\gb = \gh_3 \oplusrhrim_{ad} \ga$, where $\ga$ is the Cartan subalgebra generated by

$$U = \left(
\begin{array}{ccc}
 -1& 0 & 0\\
0 & 0 & 0\\
0 & 0& 1
\end{array} 
\right),
\;
V= \frac{1}{3}\left(
\begin{array}{ccc}
-\im & 0 & 0\\
0 & 2 \im & 0\\
0 & 0 & -\im
\end{array}
\right).
$$
If we take the basis $\{X, Y, Z\}$ for $\gh_3$:
$$X = \frac{1}{2}\left(
\begin{array}{ccc}
0& 0 & 0\\
1 & 0 & 0\\
0 & \im & 0
\end{array} 
\right),
\;
Y= \frac{1}{2}\left(
\begin{array}{ccc}
0& 0 & 0\\
\im & 0 & 0\\
0 & 1 & 0
\end{array}
\right),
\;
Z= \frac{1}{2}\left(
\begin{array}{ccc}
0& 0 & 0\\
0 & 0 & 0\\
1 & 0 & 0
\end{array}
\right)
$$
with the bracket multiplications:
$$[X,Y]=-Z, [X,Z]=[Y,Z]=0,$$
the action of $\ga$ on $\gh_3$ is given by
$$\ad_U (X)=X, \ad_U(Y)=Y, \ad_U(Z)=2Z, \ad_V(X)=Y, \ad_V(Y)=-X, \ad_V(Z)=0.$$
Note that the action of $U$ is diagonal and the action of $V$ is skew-symmetric.

We can express $\gb = \gh_3 \oplusrhrim_{ad} \ga$ as a Lie algebra extension:
$$ 0 \rightarrow \gh_3 \xrightarrow{i} \gb \xrightarrow{r} \ga \rightarrow 0,$$
where we have $\gh_3 = [\gb, \gb]$. Let $\gl$ be a three-dimensional subalgebra of $\gb$.
We have
$$ 0 \rightarrow \gl \cap \gh_3 \xrightarrow{i} \gl \xrightarrow{r} \gc \rightarrow 0,$$
where $\gl \cap \gh_3$ is a subalgebra of $\gh_3$ of dimension $3, 2$ or $1$, and $\gc$ is a
subalgebra of $\ga$ with ${\dim}\; \gc =0,1$ or $2$ accordingly.

For the case where $\gl$ is a modification of $\gh_3$ by a Lie algebra morphism $\tau: \gh_3 \rightarrow \ga$,
we have 
$$\ker \tau = \gl \cap \gh_3,\; {\rm im}\; \tau = \gc,\; \dim \ker \tau + \dim {\rm im}\; \tau = \dim \gl$$
Since $Z \in \ker \tau$, we have $\dim \gc= 1$ or $2$. In case $\dim \gc=1$, we can assume $Y \in \ker \tau$.
Then, we have $[\gl, \gl]= \gl \cap \gh_3 \simeq \bR^2$, and $\gl = [\gl,\gl] \oplusrhrim_{ad}  \bR \simeq \bR^2 \oplusrhrim_{ad}  \bR$,
where the action is diagonal.
In case $\dim \gc=2$, we have $\dim\ker \tau =1$ and thus $[\gl,\gl] =<Z>$. Since the action of $V$ on $Z$ is trivial,
we see that  $\gl \simeq \gaff(\bR) \oplus \bR$. 

For the case where $\gl$ is {\em not} any modification of $\gh_3$, we can see that $\gl$ must be either
$\gl_0=\gh_3$ or $\gl_1$;
having a basis $\{X, Z, U \}$ (or $\{Y, Z, U\}$) with bracket multiplication
$$[Z, U]=2Z, [X,U]=X, [X,Z]=0\; ({\rm or}\; [Z, U]=2Z, [Y,U]=Y, [Y,Z]=0);$$
or $\gl_2$; having a basis $\{Z, U, V\}$ with bracket multiplication
$$[Z,U]=2Z, [Z,V]=0, [U,V]=0.$$
We see that $\gl_1$ is isomorphic to $\bR^2 \oplusrhrim_{ad} \bR$, and $\gl_2$ is isomorphic to $\gaff(\bR) \oplus \bR$.

It is not hard to check that $\bR^2 \oplusrhrim_{ad}  \bR$ does not admit a non-degenerate CR structure. In fact, let $\phi$
be a contact form and $\eta$ the Reeb field with $i(\eta) \phi = 1$. We can set $\phi=p u + b y + c z$ and
$\eta = p U + b Y + c Z$, where $u, y, z$ are the Maurer-Cartan forms corresponding to 
$U, Y, Z$ respectively, and  $p, b, c \in \bR$ with $ bc \not=0, p^2+b^2+c^2=1$. Then we have $i(\eta) d \phi (U)\not= 0$.
\end{proof}


In order to determine all subalgebras $\gl$ of $\gg = \gs\gu(p+1, q+1)$  with non-degenerate CR structures
for general dimension, satisfying (6.1), note first that a non-degenerate CR structure in general 
is locally isomorphic to a contact structure, and thus for a Lie groups (algebra), it must admits an invariant contact structure.
According to a well known theorem of Boothby and Wang \cite{BW}, the only semisimple Lie algebras admitting contact structures
are $\gsu(2)$  and $\gsl(2, \bR)$. Therefore, for the determination of Cartan flat non-degenerate CR Lie algebras, it is sufficient to
consider the cases of solvable Lie algebras and show that if the Levi-decomposition of the given Lie algebra has a non-trivial
radical (the maximal solvable ideal), then it must be a solvable Lie algebra (contained in the Borel subalgebra) (see Lemma~3).

We first consider a solvable Lie subalgebra $\gl$ of the Borel subalgebra $\gb_{p,q} = \gh_{2m+1} \oplusrhrim_{ad}~\ga_{p,q}$ 
satisfying (6.1). For simplicity, we write simply $\gb, \gh, \ga$ for $\gb_{p,q}, \gh_{2m+1}, \ga_{p,q}$. 


As in the three-dimensional case, we have

$$ 0 \rightarrow \gh \xrightarrow{i} \gb \xrightarrow{r} \ga \rightarrow 0, \eqno{(6.4)}$$
with $\gh = [\gb, \gb]$, and
$$ 0 \rightarrow \gl \cap \gh \xrightarrow{i} \gl \xrightarrow{r} \gc \rightarrow 0.\eqno{(6.5)}$$

We have the standard basis $\{ X_i, Y_j, Z \}$ for $\gh$:
$$X_i = \frac{1}{2}\left(
\begin{array}{ccc}
0& 0 & 0\\
e^i & 0 & 0\\
0 & \im e_i & 0
\end{array} 
\right),
\;
Y_j= \frac{1}{2}\left(
\begin{array}{ccc}
0& 0 & 0\\
\im e^j & 0 & 0\\
0 & e_j & 0
\end{array}
\right),
\;
Z= \frac{1}{2}\left(
\begin{array}{ccc}
0& 0 & 0\\
0 & 0 & 0\\
1 & 0 & 0
\end{array}
\right),
$$
with the only non-zero bracket multiplications defined by
$$[X_i,Y_j] = -\delta_{i,j} Z, i, j = 1,2,..., m,$$
where $e^i$ is the $i$-th unit column vector of $\bC^m$ and $e_i$ is the $i$-th unit row vector of $\bC^m$.
The Cartan algebra $\ga$ is decomposed  as 
$$\ga = \gt \oplus \gs, \eqno{(6.6)}$$
where 
$$ \gt = \Big\{\left(
\begin{array}{ccc}
c & 0 & 0\\
0& 0 & 0\\
0 & 0 & -c
\end{array} 
\right) \Big| c \in \bR
\Big\},\; 
\gs = \Big\{
\left(
\begin{array}{ccc}
-\im d& 0 & 0\\
0 & \frac{2 \im d}{m} I_m + D_0& 0\\
0 & 0 & - \im d
\end{array} 
\right) \Big| d \in \bR \Big\}$$
where $D_0=\im {\rm dia}(d_1,d_2,...,d_m)$ with $d_1+d_2+ \cdots +d_m=0$.
We have a basis $\{U\}$ for $\gt$ and a basis $\{V_k\}, k=1,2,...,m$ for $\gs$ such that

$$\ad_U (X_i)=X_i, \ad_U(Y_j)=Y_j, \ad_U(Z)=2Z$$
for all $i,j=1,2,..., m$, and
$$\ad_{V_k}(X_k)=Y_k, \ad_{V_k}(Y_k)=-X_k, \ad_{V_k}(X_j)=0, \ad_{V_k}(Y_j)=0$$
for all $j \not= k$.


We first consider the case where $\gl$ is a modification of $\gh$ by a Lie algebra 
morphism $\tau: \gh \rightarrow \ga$.
In case $\ga=\gt$, that is, $\tau: \gh \rightarrow \gt$,
we have $\dim \ker \tau = \dim (\gl \cap \gh) = 2m$, and $\dim {\rm im}\, \tau =1$.
We have $\ker \tau = \gl \cap \gh = [\gl,\gl] \simeq \bR \oplus \gh_{2m-1}$ on which $U$ acts diagonally:
$$\gl \simeq (\bR \oplus \gh_{2m-1}) \oplusrhrim_{ad} \bR. \eqno(6.7)$$
In case $\ga=\gs$, that is, $\tau: \gh \rightarrow \gs$,
we have $\dim \ker \tau = \dim (\gl \cap \gh) = 2p+1$, and $\dim {\rm im}\, \tau =2q$ with $p+q=m, p \ge q$.
We can assume that $\ker \tau = \gl \cap \gh = [\gl,\gl]$ is generated by $X_{1}, Y_{1},..., X_p, Y_p, Z$, and
${\rm im}\, \tau$ is generated by $\tau(X_{p+1}), \tau(Y_{p+1}),..., \tau(X_p), \tau(Y_m)$.
As shown in the paper \cite{CH}, $\gl$ is the central extension
by $\bR$ of a meta abelian Lie algebra $\gm$:
$$ 0 \rightarrow \bR \rightarrow \gl \rightarrow \gm \rightarrow 0,$$
where 
$$ 0 \rightarrow \bR^{2p} \rightarrow \gm \rightarrow \bR^{2q} \rightarrow 0,$$
and the action of $\bR^{2q}$ on $\bR^{2p}$ is given by
$$\ad_{X_{p+i}}(X_{j}) = Y_{j},\;  \ad_{X_{p+i}}(Y_{j}) = -X_{j},\; \ad_{Y_{p+i}}(X_{j}) = Y_{j}, \; \ad_{Y_{p+i}}(Y_{j})= -X_{j},$$
where $p+q=2m$, $i=1,2,...,q, j=1,2,...,p$. Thus we have in this case
$$ \gl  \simeq \overline{\gh}_{2m+1}. \eqno{(6.8)}$$
In the general case of $\ga$, that is, $\tau: \gh \rightarrow \ga$, we consider
a composition of $\tau_t: \gh \rightarrow \gt$ and $\tau_s: \gh \rightarrow \gs$. Let $\gh_s$ be a modification of $\gh$ by $\tau_s$.
Then, we can define a well defined modification $\bar{\tau_t}: \gh_s \rightarrow \gt$ with $\bar{\tau_t} ([\gh_s, \gh_s])=0$; and
the modification by $\tau: \gh \rightarrow \ga$ is obtained as the composition $\tau_t \circ \tau_s$. We thus obtain 
$$\gl \simeq (\bR \oplus \overline{\gh}_{2m-1}) \oplusrhrim_{ad} \bR, \eqno{(6.9)}$$
where $\overline{\gh}_{2m-1}$ is a modification of ${\gh}_{2m-1}$ by $\tau_s: \gh_{2m-1} \rightarrow \gs$.

Next, we consider the  case that $\gl$ is not any modification of $\gh$. Since $\gl \cap \gh$ is a subalgebra of $\gh$ with $\dim \gl \cap \gh \ge m$,
which is invariant by the action of $\gc \subset \ga$ in (6.5), we see that $\gl$ must be either $\gl_0 = \gh$ or 
$\gl_1 = (\bR \oplus \gh_{2m-1}) \oplusrhrim_{ad} \bR$;
having a basis $\{X_1, ..., X_m, Y_1, ..., Y_{m-1}, Z, U\}$ (or  $\{X_1, ..., X_{m-1}, Y_1, ..., Y_{m}, Z, U\}$ with the bracket multiplication
$$[X_i, U]= X_i, [Y_j,U]= Y_j, [Z,U]= 2Z, [X_i, Y_j]= - \delta_{i,j} Z;$$
or $\overline{\gl}  \simeq (\bR \oplus \overline{\gh}_{2m-1}) \oplusrhrim_{ad} \bR$; having a basis  
$\{X_1, ..., X_p, Y_1, ..., Y_{p-1}, Z, U, W_1, ...,W_{2q}\}$ with $p+q=m$, where $W_j \in \gs,$
and  $\{X_1,... X_{p-1}, Y_1,..., Y_{p-1}, W_1,..., W_{2q}\}$ being a basis of $\overline{\gh}_{2m-1}$.


The integrability of complex structure on the CR structure $\mathcal D$ is preserved under modification
(c.f. \cite{CH}), and thus CR structures on $\gh_{2m+1}$ and those of $\overline{\gh}_{2m+1}$
are isomorphic. 

It is not hard to check that neither $(\bR \oplus \gh_{2m-1}) \oplusrhrim_{ad} \bR$ nor
$(\bR \oplus \overline{\gh}_{2m-1}) \oplusrhrim_{ad} \bR$ admits a non-degenerate CR structure.
In fact, let $\phi$ be a contact form and $\eta$ the Reeb field with $i(\eta) \phi = 1$. 
We can set $\phi=p u + \sum_{i=1}^{m-1} a_i x_i + \sum_{j=1}^{m}  b_j y_j+ c z$ and
$\eta = p U + \sum_{i=1}^{m-1} a_i X_i + \sum_{j=1}^{m}  b_j Y_j+ c Z$,
where $u, x_i, y_j, z$ be the Maurer-Cartan forms corresponding to 
$U, X_i, Y_j, Z$ respectively, and  $p, a_i, b_j, c \in \bR$ with 
$p^2 + \sum_{i=1}^{m-1} a_i^2 + \sum_{j=1i}^{m} b_j^2 + c^2=1$ and not all $a_i, b_j, c$ are $0$.
Then we have $i(\eta) d \phi (U)\not= 0$.


\bl Let $\gl$ be a subalgebra of $\gg=\gs\gu(p+1, q+1)$, satisfying the condition (6.1):
$$ \gg= \gl \oplus \gh_0,$$
where  $\gh_0$ is the isotropy subalgebra of $\gg$. If the radical $\gr$ of $\gl$ in the Levi-decomposition
$$\gl = \gr \oplusrhrim_{ad} \gs$$
is non-trivial, then $\gl$ must be solvable, i.e. $\gl = \gr$.
\el

\begin{proof}

Recall that $\gg$ has the canonical grading:
$$\gg = \gg^{-2} \oplus \gg^{-1} \oplus \gg^{0} \oplus \gg^{1}  \oplus \gg^{2},$$
where $\gg^{-2}, \gg^{-1}, \gg^{0}, \gg^{1}, \gg^{2}$ are the subalgebras of $\gg$ consisting of the elements

$$ X_{(a,0, 0, 0, 0,0)}, X_{(0,0, 0, 0, u, 0)}, X_{(0,0, \beta, \gu, 0, 0)}, X_{(0,0, 0, 0, 0, v)}, X_{(0,b, 0, 0, 0, 0)}$$
respectively for the elements of $\gg$ expressed in (6.3):

$$ X_{(a,b, \beta, \gu, u, v)} =\left(
\begin{array}{ccc}
\beta & - \im v^* & b\\
u & \gu_{p,q} & v\\
a &  \im \overline{u}^* & - \overline{\beta}
\end{array} 
\right),
\;\; \beta \in \bC, a, b \in \bR, u, v \in \bC^m.$$
We have the grading preserving the bracket multiplication:
$$[\gg^i, \gg^j]= \gg^{i+j},$$
with $-2 \le i, j \le 2$ and $\gg^{i+j} = \{0\}$ for $|i+j| > 2$.
Note that we have
$$\ga = \gg^{0}, \, \gh =\gg^{-2} \oplus \gg^{-1},\, \gb = \gg^{-2} \oplus \gg^{-1} \oplus \gg^{0},\, 
\gh_0 = \gg^{0} \oplus \gg^{1} \oplus \gg^{2}.$$

Now, suppose that $\gr \not= \{0\}$. We may assume  that $\gr$ is included in the Borel subalgebra $\gb$.
In case $\gr \subset \gg^{-2} \oplus \gg^{-1}$, since $\gr$ is an ideal of $\gl$, for any element $\overline{V} =
V + b \in \gl$ with $V \in \gh$ and $b \in \gh_0$, $[\gr, \overline{V}] \subset \gr$, and thus $b \in \ga$. Hence
we have $\gl \subset \gb$; that is, $\gl$ is solvable.
In case $\gr \not\subset \gg^{-2} \oplus \gg^{-1}$, take an element $\overline{U} = U + a \in \gr$ with $U \in \gh$ and
$a \in \ga$. For any $\overline{V} = V+b \in \gl$ with $V \in \gh$ and $b \in \gh_0$, we have
$$[\overline{U}, \overline{V}] = [U, V] + [a,V] + [U, b]  + [a,b],$$
which belongs to $\gr (\subset \gb)$, since $\gr$ is an ideal of $\gl$. Here, we have
$[U, V] \in \gg^{-2},\, [a,V] \in \gg^{-2} \oplus \gg^{-1}$, and in order for $[U, b]  + [a,b]$ to be in $\gb$, $b$ must be
in $\ga$, that is,  $\overline{V} \in \gb$. Hence we have $\gl \subset \gb$.
\end{proof}


We have thus shown the main theorem.

\bt
The only Lie algebras which admit Cartan flat non-degenerate CR 
structures are $\gs \gu(2), \gsl (2,\bR), \gaff (\bR) \oplus \bR$ and $\gh_{2m+1}$
with its modifications $\overline{\gh}_{2m+1}$.
\et


\noindent{\bfseries Acknowledgements}
The authors would like to express their gratitude to the referee for very careful reading of
the manuscript with some useful comments and remarks.
\\

\noindent{\bfseries Conflect of interest}
The authors have no competing interests to declare that are relevant to the content of this article.

\end{document}